\documentclass[11pt,reqno,a4paper]{amsart}

\usepackage[mode=buildmissing]{standalone}

\usepackage{amsthm,amsmath,amsfonts,amssymb,epsfig,graphicx,latexsym,float,color}
\usepackage{mathtools,multicol} 

\usepackage{enumerate,relsize,mathrsfs}
\usepackage{sidecap}
\usepackage{multirow}

\newtheorem{definition}{Definition}[section]

\newtheorem{remark}[definition]{Remark}

\newtheorem{notation}[definition]{Notation}
\newtheorem{algorithm}[definition]{Algorithm}

\usepackage{bbold}
\usepackage{bm}

\newcommand{\dx}[1]{\ensuremath{\,\mathrm{d}{#1}}}

\newcommand{\nn}{\ensuremath{\mathbb{N}}}
\newcommand{\rr}{\ensuremath{\mathbb{R}}}

\renewcommand{\vec}[1]{\ensuremath{\bm{#1}}}
\newcommand{\tens}[1]{\ensuremath{\mathcal{#1}}}
\newcommand{\mat}[1]{\ensuremath{\mathbb{#1}}}

\newcommand{\transpose}{\ensuremath{\mathsf{T}}}
\DeclareMathOperator{\abs}{abs}

\usepackage{tikz,pgfplots,pgfplotstable}
\usetikzlibrary{calc,matrix,backgrounds}
\usepgfplotslibrary{groupplots}

\pgfplotsset{compat=newest}

\pgfplotsset{
colormap={bluered}{
rgb255(0cm)=(0,0,180); rgb255(1cm)=(0,255,255); rgb255(2cm)=(100,255,0);
rgb255(3cm)=(255,255,0); rgb255(4cm)=(255,0,0); rgb255(5cm)=(128,0,0)}
}

\newcommand{\mdef}{\ensuremath{\mathbin{:=}}}

\title[Nonlinear Reduction using EGFEM]{Nonlinear Reduction using the Extended Group Finite Element Method}
\author[K.~Tolle]{Kevin Tolle$^{\dag,1}$}
\author[N.~Marheineke]{Nicole Marheineke$^1$}
\date{\today\\%
$^\dag$ \textit{Corresponding author}, email: tolle@uni-trier.de, phone: +49\,651\,201\,3472\\%
$^1$ Universit\"at Trier, FB IV - Mathematik, Lehrstuhl Modellierung und Numerik, D-54286 Trier, Germany%
}

\begin{document}

\begin{abstract}
In this paper, we develop a nonlinear reduction framework based on our recently introduced extended group finite element method. By interpolating nonlinearities onto approximation spaces defined with the help of finite elements, the extended group finite element formulation achieves a noticeable reduction in the computational overhead associated with nonlinear finite element problems. However, the problem's size still leads to long solution times in most applications. Aiming to make real-time and/or many-query applications viable, we apply model order reduction and complexity reduction techniques in order to reduce the problem size and efficiently handle the reduced nonlinear terms, respectively. For this work, we focus on the proper orthogonal decomposition and discrete empirical interpolation methods. While similar approaches based on the group finite element method only focus on semilinear problems, our proposed framework is also compatible with quasilinear problems. Compared to existing methods, our reduced models prove to be superior in many different aspects as demonstrated in three numerical benchmark problems. 
\end{abstract}

\maketitle

\noindent
\textsc{AMS-Classification:} 35G20; 65D15; 65N30; 65K99 
 

\section{Introduction} \label{sec:introduction}

Many nonlinear (finite element) problems are numerically expensive to solve with respect to computing time and storage requirements. In scenarios where many solutions need to be computed, such as for various parameter choices or initial conditions, the computational overhead quickly becomes a limiting factor. Therefore, efficient approaches that decrease the costs of solving the problem are vital. The group finite element method (GFEM) interpolates nonlinearities onto the finite element space used for the solution, which reduces overhead by allowing the time-consuming numerical integration to be precomputed and evaluating the nonlinear functions only at the nodes, see, for example, \cite{Christie1981,Fletcher1983}. Recently, we extended the GFEM in \cite{Tolle2020} to allow for better approximations and handle nonlinearities that cannot be directly interpolated onto the solution space. In particular, the extended group finite element method (EGFEM) allows for discontinuous nonlinearities, e.g., functions that depend on derivatives of the solution. It also provides more control in the approximation of the nonlinear functions, since general finite element spaces can be used instead of reusing the solution space. Additionally, we developed a tensor structure for quasilinear problems using EGFEM so that such problems can also be solved more efficiently than with the standard Galerkin approach (SGA). However, while the GFEM and EGFEM reduce the computational overhead, the size of the discrete problems is still large and impacts the solution time. Thus, projection-based model reduction techniques are crucial in order to further speed up the computations with the help of reduced order models (ROM), see, for example, \cite{Baur2014,Benner2015,Quarteroni2015,Benner2017} for more details on a variety of reduction approaches. After interpolating nonlinear functions onto finite element spaces via EGFEM, we employ the proper orthogonal decomposition (POD) method to reduce the problem's dimension with the help of a reduced basis for the solution. Finally, we apply complexity reduction through the discrete empirical interpolation method (DEIM) in order to further reduce the computational costs associated with evaluating the nonlinear terms. Our proposed approach takes advantage of the benefits of the EGFEM in comparison to the GFEM. These advantages deliver a framework that is easy to implement, directly applicable to quasilinear problems and more versatile than existing techniques. For example, the group proper orthogonal decomposition (GPOD) method presented in \cite{Dickinson2010} and the finite element method with interpolated coefficients (FEIC) for ROMs obtained from POD from \cite{Wang2015} are based on the GFEM. As such, these methods only focus on semilinear problems. Our approach -- in contrast -- is able to replicate the methods from \cite{Dickinson2010} and \cite{Wang2015}, while also being able to handle problems outside their scope. The ROMs and cROMs attained using our approach prove to be powerful approximations of the full order models (FOM), while requiring only a fraction of the computational overhead.

The structure of this paper is as follows. In Section~\ref{sec:setting} we briefly present the problem setting and introduce the notation and conventions used in Section~\ref{sec:approach}. Our three-step approach, which combines the EGFEM, POD and DEIM in order to attain efficient, reduced models, is described in Section~\ref{sec:approach}. Finally, Section~\ref{sec:results} illustrates the advantages of our method in comparison to established reduction techniques with the help of three benchmark problems.

\section{Problem Setting} \label{sec:setting}

For this work, we develop a reduction scheme for nonlinear finite element problems. Let $ \Omega \subset \rr^{d} $, $ d \in \nn $, denote a polygonal domain with boundary $ \Gamma $. Consider the following nonlinear problem:
\begin{gather} \label{eqn:model}
\begin{aligned}
    -\nabla \cdot (a(x,u(x),\nabla u(x);\mu) \, \nabla u(x)) + c(x,u(x),\nabla u(x);\mu) &= q(x;\mu), & x \in \Omega, \\
    u(x) &= 0, & x \in \Gamma,
\end{aligned}
\end{gather}
where $ \mu \in \mathcal{P} $ denotes some parameter setting. Since our focus is on the numerical solution of the finite element problem based on \eqref{eqn:model}, we assume that the scalar-valued functions $ a $, $ c $ and $ q $ are sufficiently regular and ensure the existence of a weak solution. Neumann- and Robin-type boundary conditions result in extra terms in the weak formulation but otherwise present no additional difficulties. Inhomogeneous Dirichlet conditions require purely technical adaptations to the problem formulation, while the ideas remain the same. Note that time-dependent problems can be handled analogously as seen in one of the examples considered in Section \ref{sec:results}. The finite element method introduces a finite-dimensional approximation space $ V_{h} $ using a triangulation of $ \Omega $. Assume that a nodal basis $ \{\phi_{1},\ldots,\phi_{N_{u}}\} $ exists with $ \phi_{i}(x_{j}) = \delta_{ij} $ for all degrees of freedom $ x_{j} \in \Omega $. We are interested in efficiently computing an approximate solution
\begin{gather*}
    u_{h}(x) = \sum_{i=1}^{N_{u}} u_{i} \, \phi_{i}(x),
\end{gather*}
where the coefficients $ \vec{u} = (u_{1},\ldots,u_{N_{u}})^{\transpose} $ solve the weak formulation of \eqref{eqn:model} on $ V_{h} $:
\begin{gather} \label{eqn:discrete_model}
    \mat{K}(a,\vec{u};\mu) \, \vec{u} + \vec{\ell}(c,\vec{u};\mu) = \vec{q}(\mu)
\end{gather}
with
\begin{gather} \label{eqn:forms}
\begin{aligned}
	\left[\mat{K}(a,\vec{u};\mu)\right]_{i,j} &\mdef \int_{\Omega} a(x,u_{h}(x),\nabla u_{h}(x);\mu) \, \nabla \phi_{j}(x) \cdot \nabla \phi_{i}(x) \dx{x}, \\
	\left[\vec{\ell}(c,\vec{u};\mu)\right]_{i} &\mdef \int_{\Omega} c(x,u_{h}(x),\nabla u_{h}(x);\mu) \, \phi_{i}(x) \dx{x}, \\
	\left[\vec{q}(\mu)\right]_{i} &\mdef \int_{\Omega} q(x;\mu) \, \phi_{i}(x) \dx{x}.
\end{aligned}
\end{gather}
Throughout this paper, we refer to the direct solution of \eqref{eqn:discrete_model} as the standard Galerkin approach (SGA). Many practical applications result in a large number of variables, i.e., $ N_{u} \gg 1 $. At the same time, the discrete nonlinear problem \eqref{eqn:discrete_model} is solved numerically through iterative methods, which require the repeated integration of the terms in \eqref{eqn:forms}. This is often time-intensive because of the computational overhead associated with the repeated assembly, i.e., numerical integration in \eqref{eqn:forms} for each iteration. In order to mitigate the computational costs, we investigate a combination of reduction techniques along with the EGFEM in order to derive efficient reduced models.

Before describing our reduction approach in detail, we briefly introduce the notation that is used throughout Section~\ref{sec:approach} and Section~\ref{sec:results}:
\begin{notation}[Tensor Calculus]
Vectors, matrices and third-order tensors are denoted as 
\begin{align*}
    \vec{v} &\in \rr^{N_{1}}, & \mat{M} &\in \rr^{N_{1} \times N_{2}}, & \tens{T} &\in \rr^{N_{1} \times N_{2} \times N_{3}}.
\end{align*}
For the manipulation of the third-order tensors, we use the following rules for the single and double contractions:
\begin{align*}
    \left[\tens{T} \cdot \vec{v}\right]_{i,j} &= \sum_{k} \left[\tens{T}\right]_{i,j,k} \left[\vec{v}\right]_{k}, & 
    \left[\tens{T} : \mat{M}\right]_{i} &= \sum_{j} \sum_{k} \left[\tens{T}\right]_{i,j,k} \left[\mat{M}\right]_{j,k},
\end{align*}
where the single contraction eliminates one degree of freedom and results in a matrix, while the double contraction reduces the order of the tensor by two and results in a vector. In contrast, the tensor product $ \otimes $ combines two vectors into a matrix through
\begin{gather*}
    \left[\vec{v} \otimes \vec{w}\right]_{i,j} = \left[\vec{v}\right]_{i} \, \left[\vec{w}\right]_{j}. 
\end{gather*}
For the reduction, we also use the $ n $-mode product $ \times_{n} $ for $ n \in \{1,2,3\} $: 
\begin{align*}
    \left[\tens{T} \times_{1} \mat{A}_{1}\right]_{\tilde{i},j,k} &= \sum_{i} \left[\tens{T}\right]_{i,j,k} \, \left[\mat{A}_{1}\right]_{\tilde{i},i}, & \mat{A}_{1} &\in \rr^{\tilde{N}_{1} \times N_{1}}, \\
    \left[\tens{T} \times_{2} \mat{A}_{2}\right]_{i,\tilde{j},k} &= \sum_{j} \left[\tens{T}\right]_{i,j,k} \, \left[\mat{A}_{2}\right]_{\tilde{j},j}, & \mat{A}_{2} &\in \rr^{\tilde{N}_{2} \times N_{2}}, \\
    \left[\tens{T} \times_{3} \mat{A}_{3}\right]_{i,j,\tilde{k}} &= \sum_{k} \left[\tens{T}\right]_{i,j,k} \, \left[\mat{A}_{3}\right]_{\tilde{k},k}, & \mat{A}_{3} &\in \rr^{\tilde{N}_{3} \times N_{3}},
\end{align*}
where the result is a third-order tensor with the $ n $-th dimension changed according to the matrix used in the product. The definition of the $ n $-mode product is based on \cite{Kolda2009}, which we refer the reader to for more details.
\end{notation}

\section{Reduction Approach} \label{sec:approach}

Model order reduction techniques serve a vital role in efficient solution strategies. Many nonlinear reduction approaches rely on data-driven methods. However, it has been shown that model order reduction alone is not effective, since nonlinear functions still depend on the full problem dimension, see, for example, \cite{Chaturantabut2010} for details. Therefore, complexity reduction in the form of (D)EIM \cite{Barrault2004,Chaturantabut2010} or empirical cubature \cite{Hernandez2017} is necessary in order to efficiently evaluate the reduced nonlinear terms. In this paper, we focus on POD for the model reduction and DEIM for the complexity reduction. In combination with the EGFEM, our scheme delivers a powerful tool for the reduction of nonlinear finite element models that is more versatile than the approaches described in \cite{Dickinson2010,Wang2015}. For example, discontinuous nonlinearities can be handled efficiently with the help of EGFEM. When the DEIM is applied to point-wise schemes, such as finite difference discretizations, the complexity reduction is non-intrusive. This advantage is carried over for the reformulation via GFEM and EGFEM, since the nonlinearities are evaluated point-wise. In contrast, applying the DEIM to the SGA requires modifications to the finite element implementation in order to ensure that the complexity reduction actually improves the evaluation of the reduced nonlinear terms. Furthermore, the tensor structure introduced in \cite{Tolle2020} allows for a straight-forward reduction of the quasilinear terms. As such, our approach can be seen as an extension of \cite{Dickinson2010,Wang2015}, which uses the GFEM to approximate nonlinear terms before reducing with POD and optionally DEIM. The reduction approach can be separated into the following three steps, which are explained in more detail in what follows:
\begin{enumerate}
	\item apply the EGFEM in order to rewrite nonlinear terms,
	\item compute projection matrices using POD in order to derive a ROM,
	\item use the DEIM to attain a cROM, which efficiently handles reduced nonlinear terms.
\end{enumerate}

\subsection{Extended Group Finite Element Method}

Most of the computational overhead for nonlinear finite element problems lies in the assembly of nonlinear forms. Because the problems must be solved using iterative methods, such as Newton's methods, these forms have to be assembled in every iteration. Afterwards, linear systems are solved and the approximate solution is updated. In order to reduce the overhead and consequently speed up the computations, the (extended) group finite element method interpolates nonlinearities onto finite element spaces. The major difference between the GFEM and EGFEM is the choice of the approximation space for the nonlinear functions, i.e., GFEM interpolates nonlinear terms onto the solution space $ V_{h} $, while EGFEM introduces an additional finite element space $ W_{h} $ that can be tailored to the nonlinearity. As such, EGFEM encapsulates GFEM for $ W_{h} = V_{h} $ but also allows for more accurate approximations and the treatment of more general nonlinearities, such as functions that depend on derivatives of the solution. For a more detailed investigation into the EGFEM including guidelines and numerical studies, we refer the reader to \cite{Tolle2020}.

The basic idea is to introduce a finite element approximation space $ W_{h} $ with a nodal basis $ \{\eta_{1}^{f}, \ldots, \eta_{N_{f}}^{f}\} $ that is used to interpolate the nonlinear function $ f(u_{h},\nabla u_{h};\mu) $. Instead of evaluating the nonlinear terms in \eqref{eqn:discrete_model} in every iteration, we approximate 
\begin{gather*}
    f(x,u_{h}(x),\nabla u_{h}(x);\mu) \approx \sum_{j=1}^{N_{f}} f_{j}(\vec{u};\mu) \, \eta_{j}^{f}(x).
\end{gather*}
The coefficients $ f_{j}(\vec{u};\mu) $, $ j = 1, \ldots, N_{f} $, are determined by evaluating the function at the degrees of freedom $ x_{j}^{f} $, i.e., $ f_{j}(\vec{u};\mu) \mdef f(x_{j}^{f},u_{h}(x_{j}^{f}),\nabla u_{h}(x_{j}^{f});\mu) $. In the case of GFEM, we have the special choice of $ W_{h} = V_{h} $ and, therefore, $ \eta_{i}^{f} = \phi_{i} $, $ i = 1,\ldots,N_{u} $, with the coefficients evaluated at the same degrees of freedom as the solution. For the nonlinear terms in \eqref{eqn:forms}, this procedure results in the approximations
\begin{align*}
    \left[\mat{K}(a,\vec{u};\mu)\right]_{ij} &\approx \left[\tens{K}^{a} \cdot \vec{a}(\vec{u};\mu)\right]_{i,j}, & \tens{K}^{a} &\mdef \int_{\Omega} \eta_{k}^{a}(x) \, \nabla \phi_{j}(x) \cdot \nabla \phi_{i}(x) \dx{x}, \\
	\left[\vec{\ell}(c,\vec{u};\mu)\right]_{i} &\approx \left[\mat{L}^{c} \, \vec{c}(\vec{u};\mu)\right]_{i}, & \mat{L}^{c} &\mdef \int_{\Omega} \eta_{j}^{c}(x) \, \phi_{i}(x) \dx{x},
\end{align*}
where the integrals can be precomputed and result in the third-order tensor $ \tens{K}^{a} $ and matrix $ \mat{L}^{c} $, respectively. By introducing these approximations to \eqref{eqn:discrete_model}, the EGFEM solves the following problem:
\begin{gather} \label{eqn:egfem}
    \tens{K}^{a} : (\vec{u} \otimes \vec{a}(\vec{u};\mu)) + \mat{L}^{c} \, \vec{c}(\vec{u};\mu) = \vec{q}(\mu)
\end{gather}
with
\begin{align*}
    \vec{a}(\vec{u};\mu) &= A(\vec{x}^{a},\Pi^{a}_{u} \, \vec{u},\Pi^{a}_{\nabla u} \cdot \vec{u};\mu), \\
    \vec{c}(\vec{u};\mu) &= C(\vec{x}^{c},\Pi^{c}_{u} \, \vec{u},\Pi^{c}_{\nabla u} \cdot \vec{u};\mu),
\end{align*}
where the coefficient vectors $ \vec{a} $ and $ \vec{c} $ are defined point-wise. This leads to the introduction of the functions $ F : \rr^{d \times N_{f}} \times \rr^{N_{f}} \times \rr^{d \times N_{f}} \to \rr^{N_{f}} $, $ (F,f) \in \{(A,a),(C,c)\} $, which are defined through 
\begin{align*}
    (\vec{x}^{f},\Pi_{u}^{f} \, \vec{u},\Pi_{\nabla u}^{a} \cdot \vec{u};\mu) &\mapsto (f(x_{1}^{f},u_{h}(x_{1}^{f}),\nabla u_{h}(x_{1}^{f});\mu),\ldots,f(x_{N_{f}}^{f},u_{h}(x_{N_{f}}^{f}),\nabla u_{h}(x_{N_{f}}^{f});\mu))^{\transpose},
\end{align*}
as well as the interpolation operators $ \Pi_{u}^{f} : \rr^{N_{u}} \to \rr^{N_{f}} $ and $ \Pi_{\nabla u}^{f} : \rr^{N_{u}} \to \rr^{d \times N_{f}} $, which are used to evaluate the function $ u_{h} $ and its gradient $ \nabla u_{h} $ at the degrees of freedom $ \vec{x}^{f} = (x_{1}^{f},\ldots,x_{N_{f}}^{f})^{\transpose} $ associated to the approximation space $ W_{h}^{f} $ for $ f \in \{a,c\} $.
\begin{remark}
    Note that the extended group finite element approximation can also be used for nonlinear boundary terms. In this case, the approximation space only has to be defined on the boundary, which results in an additional dimension reduction. Furthermore, it can also be applied to source terms in time-dependent problems, which have to be reassembled in each time step.
\end{remark}

\subsection{Proper Orthogonal Decomposition}

The next major reduction in complexity is achieved by applying the proper orthogonal decomposition technique on the spatial dependence of the system. This means that the size of the discrete problem is reduced with the help of orthogonal projections that select the most important modes of information with the help of a set of snapshots, i.e. a collection of solutions at different times and/or different parameter settings, see, for example, \cite{Quarteroni2015}. The snapshots are collected in the matrix
\begin{align*}
	\mat{Y} &= [ \vec{u}_{1}, \ldots, \vec{u}_{n_{s}} ] \in \rr^{N_{u} \times n_{s}},
\end{align*}
where $ n_{s} $ denotes the number of snapshots and $ \vec{u}_{k} $ represents the solution in a specific configuration, e.g., for a fixed parameter value and/or at a specific time. The singular value decomposition (SVD) of $ \mat{Y} $ delivers
\begin{align*}
	\mat{Y} &= \mat{U} \, \mat{\Sigma} \, \mat{W}^{\transpose}
\end{align*}
with $ \mat{U} \in \rr^{N_{u} \times N_{u}} $ and $ \mat{W} \in \rr^{n_{s} \times n_{s}} $. Then, the projection matrix $ \mat{V}_{u} \in \rr^{N_{u} \times n_{u}} $ consists of the first $ n_{u} \ll N_{u} $ columns from the left projection $ \mat{U} $.

For our problem, the snapshot matrix consists of the coefficient vectors $ \vec{u} $ at discrete parameter settings, i.e., $ \vec{u}_{i} \mdef \vec{u}(\mu_{i}) $ for $ 1 \leq i \leq n_{s} $. After computing the SVD of $ \mat{Y} $ and determining the projection matrix $ \mat{V}_{u} $, the reduced-order variable $ \vec{u}_{\mathsf{r}} \in \rr^{n_{u}} $ with $ n_{u} \ll N_{u} $ approximates the full-order variable $ \vec{u} $ through
\begin{gather*}
	\vec{u} \approx \mat{V}_{u} \, \vec{u}_{\mathsf{r}}.
\end{gather*}
The projection-based reduction of $ \vec{u} $ in the SGA \eqref{eqn:discrete_model} and EGFEM \eqref{eqn:egfem} results in the following respective ROMs:
\begin{alignat}{4} \label{eqn:reduced_sga}
	\mat{V}_{u}^{\transpose} \, \mat{K}(a,\mat{V}_{u} \, \vec{u}_{\mathsf{r}};\mu) \, \mat{V}_{u} \, \vec{u}_{\mathsf{r}} ~&+&~ \mat{V}_{u}^{\transpose} \, \vec{\ell}(c,\mat{V}_{u} \, \vec{u}_{\mathsf{r}};\mu) ~&=&~ \vec{q}_{\mathsf{r}}(\mu), \\ \label{eqn:reduced_egfem}
	\mat{V}_{u}^{\transpose} \, \tens{K}^{a} : (\mat{V}_{u} \, \vec{u}_{\mathsf{r}} \otimes \vec{a}(\mat{V}_{u} \, \vec{u}_{\mathsf{r}};\mu)) ~&+&~ \mat{V}_{u}^{\transpose} \, \mat{L}^{c} \, \vec{c}(\mat{V}_{u} \, \vec{u}_{\mathsf{r}};\mu) ~&=&~ \vec{q}_{\mathsf{r}}(\mu)
\end{alignat}
with $ \vec{q}_{\mathsf{r}}(\mu) = \mat{V}_{u}^{\transpose} \, \vec{q}(\mu) $. Irregardless of which reduced model is considered, a major hindrance is the evaluation of the nonlinear functions. This is because these functions require the prolongation of the reduced variables to the full order dimension, evaluation of the function and finally the restriction back to the reduced order dimension. Therefore, the evaluation of the nonlinearities still depends on the size of the full order model. In light of this, we use the DEIM in order to remove the dependence on the full order dimension.

\subsection{Discrete Empirical Interpolation Method}

The difficulty with nonlinear reduction is that the evaluation of the nonlinear functions depends on the full-order dimension. In \cite{Chaturantabut2010}, the discrete empirical interpolation method solves this problem via an interpolation approach, which only uses a few components of the nonlinear function as an approximation.

In the scope of this work, we consider $ \vec{f}(\mu) \mdef \vec{f}(\mat{V}_{u} \, \vec{u}_{\mathsf{r}};\mu) $, where $ \vec{f} $ comes from either the nonlinear terms $ \mat{K}(a,\vec{u};\mu) $ and $ \vec{\ell}(c,\vec{u};\mu) $ in the case of the SGA or the point-wise evaluations $ \vec{a}(\vec{u};\mu) $ and $ \vec{c}(\vec{u};\mu) $ in the EGFEM setting for every $ \mu \in \mathcal{P} $. Note that the matrices $ \mat{K}(a,\vec{u};\mu) \in \rr^{N_{u} \times N_{u}} $ are treated as vectors in $ \rr^{N_{u}^{2}} $ and reduced as described here. In this case, one speaks of the matrix discrete empirical interpolation method (MDEIM), see \cite{Carlberg2012,Wirtz2014}. Similar to the reduction of $ \vec{u} $, the form of the approximation of $ \vec{f} $ is given through
\begin{gather*}
    \vec{f}(\tau) \approx \mat{V}_{f} \, \vec{f}_{\mathsf{r}}(\tau),
\end{gather*}
where $ \tau $ is either a parameter setting $ \mu $ or a time $ t $, with $ \mat{V}_{f} \in \rr^{N_{f} \times n_{f}} $ and $ \vec{f}_{\mathsf{r}}(\tau) \in \rr^{n_{f}} $. However, the system $ \vec{f}(\tau) = \mat{V}_{f} \, \vec{f}_{\mathsf{r}}(\tau) $ is overdetermined. Therefore, in order to get a well-defined approximation, which is efficient to evaluate, $ n_{f} $ unique rows must be selected in order to determine $ \vec{f}_{\mathsf{r}} $. This is accomplished with the help of a selection matrix
\begin{gather*}
    \mat{P}_{f} = [\vec{e}_{{p}_{1}}, \ldots, \vec{e}_{{p}_{n_{f}}}] \in \rr^{N_{f} \times n_{f}},
\end{gather*}
where $ \vec{e}_{i} $ denotes the $ i $-th canonical basis vector in $ \rr^{N_{f}} $. Then, the coefficient vector $ \vec{f}_{\mathsf{r}} $ is uniquely determined through $ \mat{P}_{f}^{\transpose} \, \vec{f}(\tau) = (\mat{P}_{f}^{\transpose} \, \mat{V}_{f}) \, \vec{f}_{\mathsf{r}}(\tau) $. This leads to the approximation
\begin{align*}
    \vec{f}(\tau) &\approx \mat{D}_{f} \, \mat{P}_{f}^{\transpose} \, \vec{f}(\tau), & \mat{D}_{f} &\mdef \mat{V}_{f} \, (\mat{P}_{f}^{\transpose} \, \mat{V}_{f})^{-1}
\end{align*}
with the DEIM operator $ \mat{D}_{f} $, where the projection matrix $ \mat{V}_{f} $ and the selection matrix $ \mat{P}_{f} $ still need to be specified. POD is applied to the nonlinear snapshots
\begin{gather*}
    \mat{Y}_{f} = [\vec{f}(\vec{u}_{1}),\ldots,\vec{f}(\vec{u}_{n_{s}})] \in \rr^{N_{f} \times n_{s}},
\end{gather*}
which results in the projection matrix $ \mat{V}_{f} $, while the selection matrix $ \mat{P}_{f} $ is chosen inductively using a Greedy-type algorithm from \cite{Chaturantabut2010}. The DEIM-algorithm based on our notation is presented in Algorithm~\ref{alg:deim}, where the first two steps compute the projection matrix $ \mat{V}_{f} $ via POD and the rest of the algorithm selects the rows for $ \mat{P}_{f} $ via a Greedy algorithm.
\begin{algorithm}[DEIM] \label{alg:deim}
    Let $ \mat{Y}_{f} $ and $ 1 \leq n_{f} \leq \min\{n_{s},N_{f}\} $ be given. 
    \begin{enumerate}[1.]
        \item Compute the SVD of $ \mat{Y}_{f} = \mat{U}_{f} \, \mat{\Sigma}_{f} \, \mat{W}_{f}^{\transpose} $.
        \item Set $ \mat{V}_{f} = [\vec{u}_{f}^{1},\ldots,\vec{u}_{f}^{n_{f}}] $, where $ \vec{u}_{f}^{i} $ denotes the $ i $-th column of $ \mat{U}_{f} $.
        \item Determine $ p_{1} = {\arg\max}_{1 \leq i \leq N_{f}} \abs([\vec{u}_{f}^{1}]_{i}) $. Set $ \mat{P}_{f} = [\vec{e}_{p_{1}}] $.
        \item Iterate over the remaining vectors $ \vec{u}_{f}^{\ell} $ with $ 2 \leq \ell \leq n_{f} $:
        \begin{enumerate}[i)]
            \item Solve $ \mat{P}_{f}^{\transpose} \, \vec{u}_{f}^{\ell} = (\mat{P}_{f}^{\transpose} \, \mat{V}_{f,\ell-1}) \, \vec{u}_{f,\mathsf{r}}^{\ell} $ for $ \vec{u}_{f,\mathsf{r}}^{\ell} $ with $ \mat{V}_{f,\ell-1} = [\vec{u}_{f}^{1},\ldots,\vec{u}_{f}^{\ell-1}] $.
            \item Determine $ p_{\ell} = {\arg\max}_{1 \leq i \leq N_{f}} \abs(\left[\vec{u}_{f}^{\ell} - \mat{V}_{f,\ell-1} \, \vec{u}_{f,\mathsf{r}}^{\ell}\right]_{i}) $. Set $ \mat{P}_{f} = [\mat{P}_{f},\vec{e}_{p_{\ell}}] $.
        \end{enumerate}
    \end{enumerate}
\end{algorithm}

For our model problem, the DEIM is applied to \eqref{eqn:reduced_sga} and \eqref{eqn:reduced_egfem} in order to efficiently evaluate the nonlinear terms. Beginning with the SGA-ROM in \eqref{eqn:reduced_sga}, DEIM approximates the nonlinear functions $ \mat{K}(a,\vec{u};\mu) $ and $ \vec{\ell}(c,\vec{u};\mu) $. This leads to the following SGA-cROM:
\begin{gather} \label{eqn:crom_sga}
    \mat{V}_{u}^{\transpose} \, T^{-1}(\mat{D}_{\mat{K}} \, \mat{P}_{\mat{K}}^{\transpose} \, T(\mat{K}(a,\mat{V}_{u} \, \vec{u}_{\mathsf{r}};\mu))) \, \mat{V}_{u} \, \vec{u}_{\mathsf{r}} + \mat{V}_{u}^{\transpose} \, \mat{D}_{\vec{\ell}} \, \mat{P}_{\vec{\ell}}^{\transpose} \, \vec{\ell}(c,\mat{V}_{u} \, \vec{u}_{\mathsf{r}};\mu) = \vec{q}_{\mathsf{r}}(\mu),
\end{gather}
where the functions $ T $ and $ T^{-1} $ denote the transformation of a matrix to a vector (column-wise) and vice versa. For the EGFEM-ROM in \eqref{eqn:reduced_egfem}, we attain the EGFEM-cROM:
\begin{gather} \label{eqn:crom_egfem}
    \tens{K}^{a}_{\mathsf{r}} : (\vec{u}_{\mathsf{r}} \otimes \mat{P}_{\vec{a}}^{\transpose} \, \vec{a}(\mat{V}_{u} \, \vec{u}_{\mathsf{r}};\mu)) + \mat{V}_{u}^{\transpose} \, \mat{L}^{c} \, \mat{D}_{\vec{c}} \, \mat{P}_{\vec{c}}^{\transpose} \, \vec{c}(\mat{V}_{u} \, \vec{u}_{\mathsf{r}};\mu) = \vec{q}_{\mathsf{r}}(\mu)
\end{gather}
with
\begin{align*}
    \tens{K}^{a}_{\mathsf{r}} &\mdef \tens{K}^{a} \times_{1} \mat{V}_{u} \times_{2} \mat{V}_{u} \times_{3} \mat{D}_{\vec{a}}.
\end{align*}
The remarkable difference between \eqref{eqn:crom_sga} and \eqref{eqn:crom_egfem} is the ``hidden'' overhead contained in \eqref{eqn:crom_sga}. While the numerical integration in \eqref{eqn:crom_egfem} is only performed once in the computation of $ \tens{K}^{a} $ and $ \mat{L}^{c} $, the matrix $ \mat{K}(a,\mat{V}_{u} \vec{u}_{\mathsf{r}};\mu) $ and vector $ \vec{\ell}(c,\mat{V}_{u} \, \vec{u}_{\mathsf{r}};\mu) $ require assembly for each evaluation. In light of this, an effective implementation of \eqref{eqn:crom_sga} requires a modification of the assembly process in order to take full advantage of the complexity reduction from the DEIM, so that information that is immediately disregarded by the selection operator $ \mat{P}^{\transpose}_{f} $, $ f \in \{\mat{K},\vec{\ell}\} $, is not generated at all. In contrast, the evaluations of $ \vec{a} $ and $ \vec{c} $ in \eqref{eqn:crom_egfem}, which are point-wise, can take full advantage of the DEIM. Similar to the DEIM for finite difference methods in \cite{Chaturantabut2010}, the selection operator $ \mat{P}_{f}^{\transpose} $, $ f \in \{\vec{a},\vec{c}\} $, can be applied directly to the arguments of $ f $. This means that the nonlinear functions are only evaluated at the $ n_{f} $ points chosen in the DEIM algorithm.

\section{Numerical Results} \label{sec:results}

In this section, we focus on the quality of the (complexity) reduced models in comparison with the FOMs. We consider three examples that focus on different aspects of the considered method(s). The first example is taken from \cite{Chaturantabut2010} and compares the SGA, GFEM and EGFEM for a semilinear problem. For the second example, we solve the viscous Burgers' equation in two dimensions, which is also considered in \cite{Dickinson2010}. This illustrates how our approach can also be applied to time-dependent problems. Finally, the third example highlights the superiority of our approach in comparison to the existing methods based on the GFEM. Note that all examples assume homogeneous Dirichlet boundary conditions. For all the examples, we use piecewise polynomial functions for the approximation space, where $ P_{k} $ denotes the piecewise polynomial functions of degree $ k \in \nn_{0} $. The solution space $ V_{h} $ uses piecewise linear functions, i.e., $ P_{1} $. The interpolation of the nonlinear function onto linear finite elements corresponds to the GFEM, since $ V_{h} = W_{h} $. In contrast, we consider the approximation spaces $ W_{h} $ using piecewise constant and quadratic functions, which are denoted with $ P_{0} $ and $ P_{2} $, respectively, for the EGFEM. The ROMs are attained using POD, while the cROMs use POD in conjunction with the DEIM. The size of the complexity reduction, i.e., the number of DEIM modes retained, is provided in parenthesis, where applicable. Finally, all computations are performed in Matlab (R2020a) on a server with a Xeon E5-2699V4 processor and 756 GB of memory. The triangulations are generated using the software Gmsh \cite{Geu09}, while an external toolbox for tensors \cite{TTB_Sparse,TTB_Software} is used for sparse third-order tensors. The stationary examples are solved using Matlab's \texttt{fsolve} function with all tolerances set to $ 10^{-8} $ and using the initial guess $ \vec{u}_{0} = \vec{0} $. The instationary (second) example is solved using Matlab's \texttt{ode15s} method with backward differentiation formula, where the absolute and relative tolerances are set to $ 10^{-8} $ and $ 10^{-6} $, respectively. In order to ensure optimal performance of Matlab's solvers, we provide the Jacobian matrix analytically. Finally, the runtimes refer to the average time elapsed as measured by Matlab's \texttt{tic}-\texttt{toc} functionality.

\subsection{Parameter-dependent Semilinear Equation}

The first example is a nonlinear parameter-dependent problem with a homogeneous Dirichlet boundary condition from \cite{Chaturantabut2010}. The equation considered is given by
\begin{align*}
    -\nabla \cdot (\nabla u(x)) + \frac{\mu_{1}}{\mu_{2}} \left(e^{\mu_{2} \, u(x)} - 1\right) &= 100 \, \sin(2 \, \pi \, x_{1}) \, \sin(2 \, \pi \, x_{2}), & x = (x_{1},x_{2}) \in [0,1]^{2}
\end{align*}
with the parameter setting $ \mu = (\mu_{1},\mu_{2})^{\transpose} \in \mathcal{P} = [0.01,10]^{2} $. With respect to \eqref{eqn:model}, we have
\begin{align*}
    a(x) &\equiv 1, & 
    c(u;\mu) &= \frac{\mu_{1}}{\mu_{2}} \left(e^{\mu_{2} \, u} - 1\right), & 
    q(x,\mu) &= 100 \, \sin(2 \, \pi \, x_{1}) \, \sin(2 \, \pi \, x_{2}).
\end{align*}
and the FOMs as well as ROMs and cROMs are given through \eqref{eqn:discrete_model}, \eqref{eqn:reduced_sga} and \eqref{eqn:crom_egfem} for the SGA and \eqref{eqn:egfem}, \eqref{eqn:reduced_egfem} and \eqref{eqn:crom_egfem} for the EGFEM. Finally, the derivative of the nonlinearity with respect to $ \vec{u} $ for the SGA and EGFEM is given through
\begin{align*}
    \left[\mathrm{D}_{\vec{u}} \vec{\ell}(c,\vec{u};\mu)\right]_{i,j} &= \int_{\Omega} \mu_{1} \, \exp\!\left(\mu_{2} \, u_{h}(x)\right) \phi_{j}(x) \, \phi_{i}(x) \dx{x}, \\
    \mathrm{D}_{\vec{u}} \vec{c}(\vec{u};\mu_{1},\mu_{2}) &= \mu_{1} \, \mathrm{diag}(\exp(\mu_{2} \, \Pi_{u}^{c} \, \vec{u})) \, \Pi_{u}^{c},
\end{align*}
where the exponential function acts on each element of its vector-valued argument and $ \mathop{\mathrm{diag}}(\cdot) $ denotes the matrix with the argument on its diagonal. The reduced models have similar derivatives with the inclusion of the chain rule, which results in additional scaling with the projection matrix $ \mat{V}_{u} $.

\begin{remark}
Note that the constant factor in $ c $ could also be shifted to the source term $ q $. However, the current choice matches with \cite{Chaturantabut2010} and ensures that $ c \equiv 0 $ on the boundary.
\end{remark}

\begin{figure}[t] \centering
\includegraphics[width=0.8\linewidth]{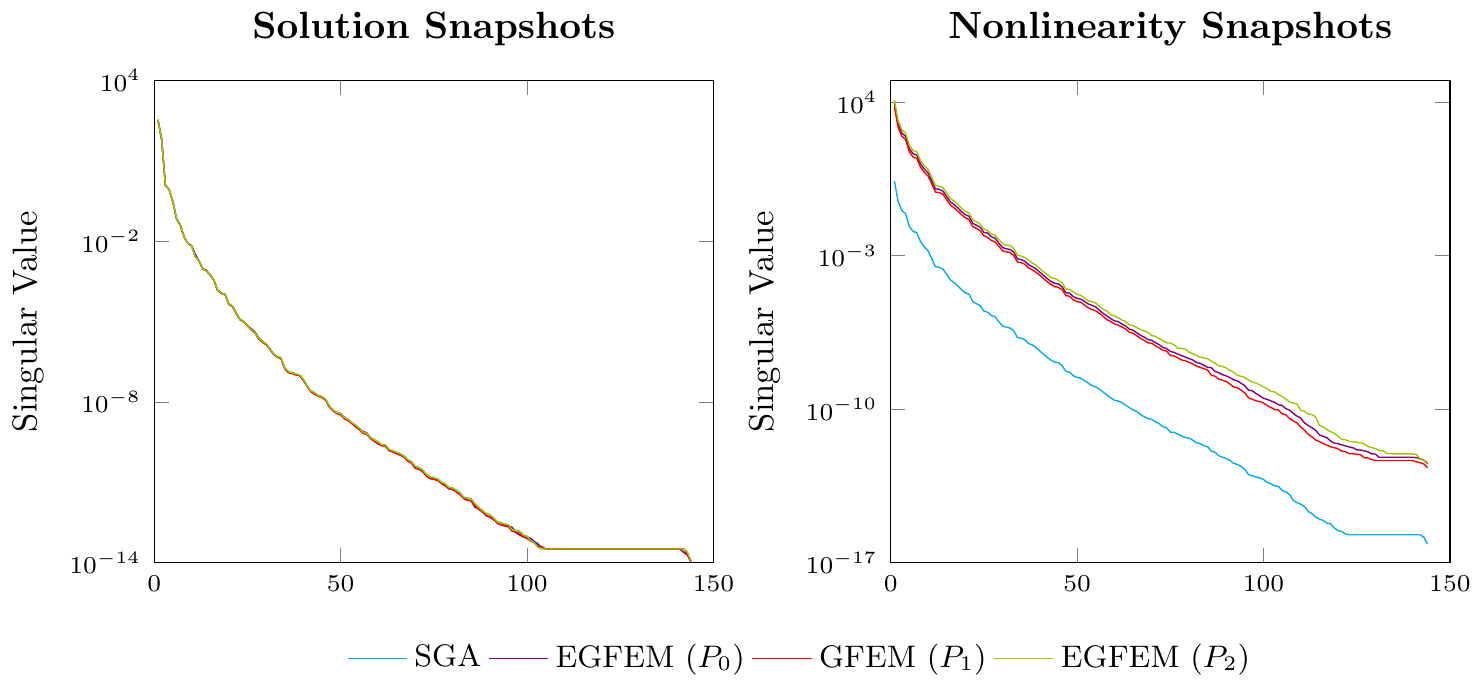}
\caption{Singular values associated to the snapshot matrices for the solution (\textsl{left}) and the nonlinearity (\textsl{right}).}
\label{fig:parameterModes}
\end{figure}

\begin{figure} \centering
\includegraphics[width=0.75\linewidth]{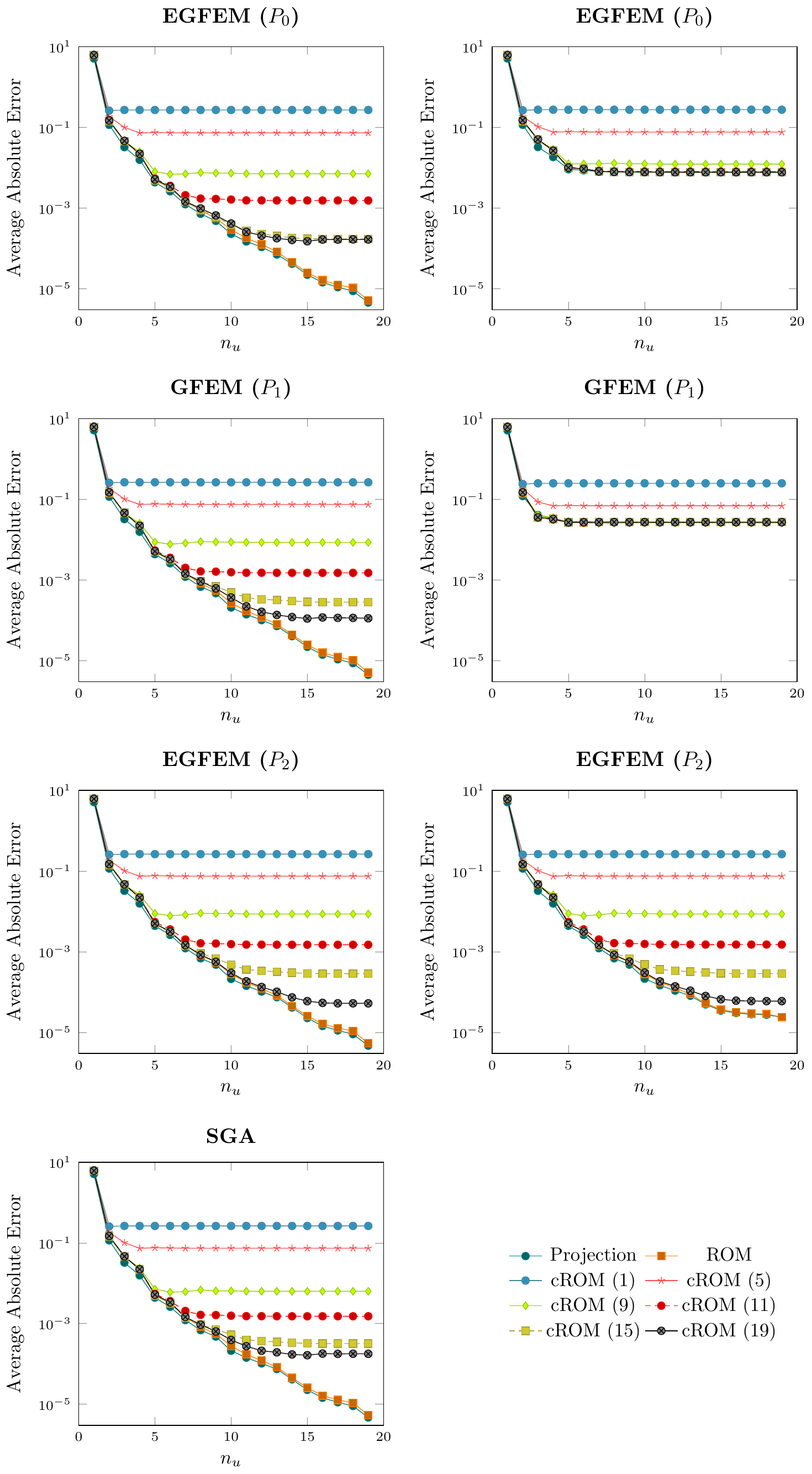}
\caption{Average absolute error over all 225 parameter samples with respect to the FOM (\textsl{left}) and to the SGA-FOM (\textsl{right}).}
\label{fig:parameterError}
\end{figure}

\begin{figure} \centering
\includegraphics[width=0.8\linewidth]{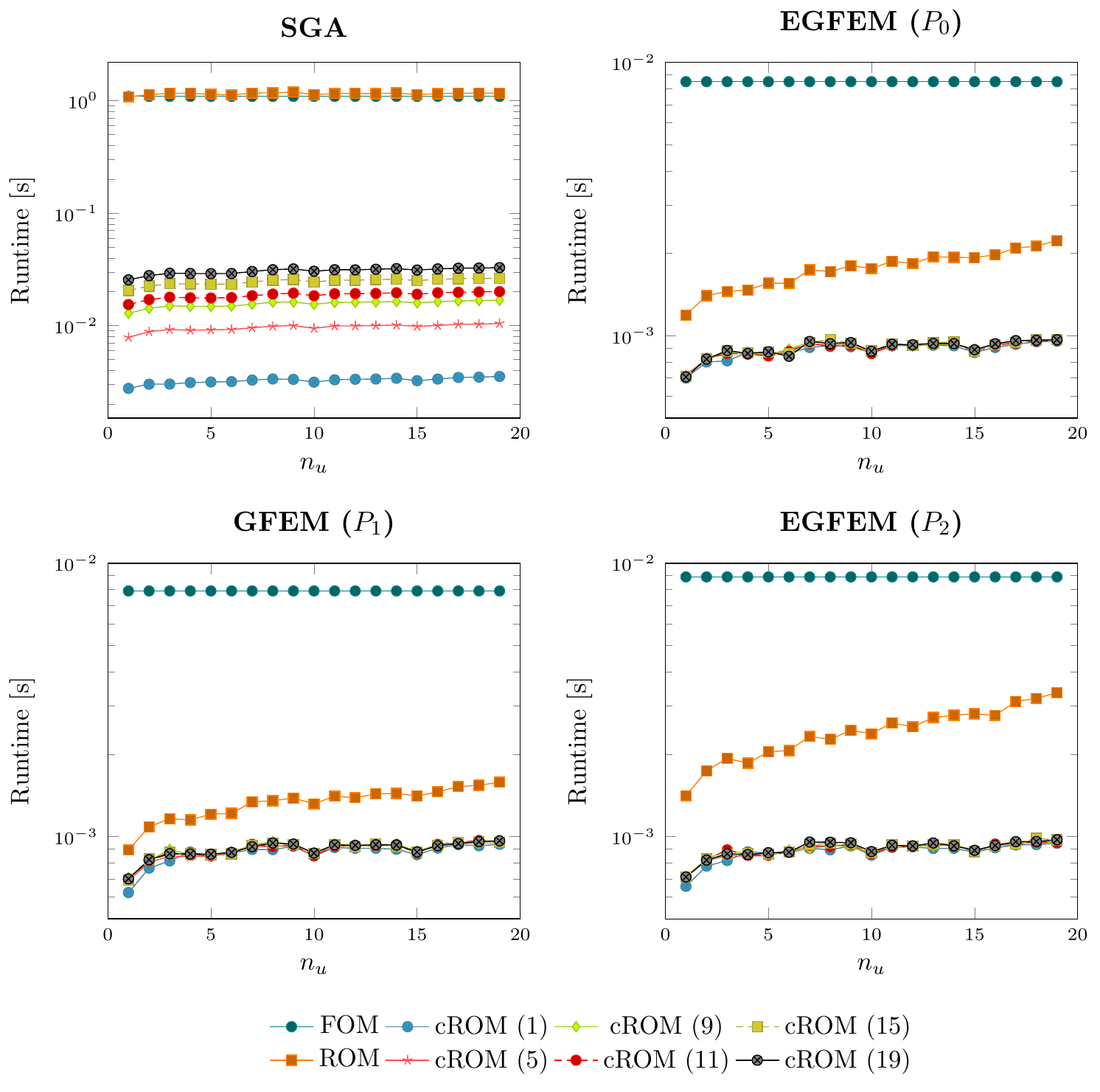}
\caption{Average time to solve the respective system.}
\label{fig:parameterTimes}
\end{figure}

Following \cite{Chaturantabut2010}, the solution and nonlinearity snapshots are generated by sampling the parameter space $ \mathcal{P} = [0.01,10]^{2} $ uniformly with 144 samples. The singular values belonging to the snapshot matrices are shown in Figure~\ref{fig:parameterModes}. Up to a constant shift, the singular values of the nonlinear snapshot matrices show almost the same decay rate. This shift in the singular values of the nonlinear snapshots between the SGA and EGFEM problems is most likely due to the additional integration step involved in the SGA. We evaluate the ROMs and cROMs on a uniform sampling of the parameter space using 225 samples, which do not coincide with the training parameters. In order to compare the different choices of $ W_{h} $ in the EGFEM and because we do not have an analytical solution available, we consider the average absolute error using the Euclidean norm with respect to two different references. The average absolute errors of the reduced solutions with respect to the full order solution are shown in the left column of Figure~\ref{fig:parameterError}. These errors are comparable to those shown in \cite{Chaturantabut2010}. The right column in Figure~\ref{fig:parameterError} shows the error with respect to the SGA-FOM. This choice as a reference is motivated by the fact that the SGA coupled with a fourth-order numerical quadrature rule for the nonlinear term should be the most accurate approximation of the exponential function in the nonlinearity. In contrast, the GFEM and EGFEM use lower-degree polynomial approximations of the nonlinear function $ c $ in order to simplify the problem. From the right column of Figure~\ref{fig:parameterError}, we see that the GFEM ($ P_{1} $) approximation as well as the EGFEM approximation with piece-wise constant functions ($ P_{0} $) poorly describe the reference SGA solution, while the EGFEM with piece-wise quadratic functions ($ P_{2} $) captures the behavior of the SGA solution and even performs better in some configurations, e.g., when using 19 DEIM modes. Finally, we see noticeable reductions in the computational overhead for the EGFEM formulations in comparison to the SGA, as shown in Figure~\ref{fig:parameterTimes}. The SGA-ROMs require the complexity reduction via DEIM (SGA-cROMs) in order to even show an improvement in the computing times. Our approach with the EGFEM -- in contrast to the SGA -- turns out to be even faster than most of the SGA-cROMs. Additionally, DEIM reduces the average time per iteration by an additional order of magnitude for the models using the EGFEM. Interestingly, as seen in Figure~\ref{fig:parameterTimes}, the average time per iteration for the GFEM- and EGFEM-cROMs appears to be independent of the DEIM dimension. This is most likely due to the fact that Matlab efficiently handles matrix-vector multiplications, and the difference between the DEIM dimensions is negligible at the considered scale.

\subsection{Burgers' Equation}

For the second example, we consider the viscous Burgers' equation in two dimensions with a homogeneous Dirichlet boundary condition:
\begin{gather*}
    \partial_{t} u(x,t) - \nabla \cdot (\nu \, \nabla u(x,t)) + \frac{1}{2} \left( \partial_{x_{1}} u(x,t)^{2} + \partial_{x_{2}} u(x,t)^{2} \right) = q(x,t)
\end{gather*}
on the domain $ x = (x_{1},x_{2}) \in [0,1]^{2} $ and time interval $ t \in [0,10] $ with $ u(x,0) = u_{0}(x) $ and a constant viscosity $ \nu = 1/100 $. In the spirit of \cite{Dickinson2010}, we choose the source $ q $ via the method of manufactured solutions and the initial condition $ u_{0} $, so that the solution is given through
\begin{gather*}
    \hat{u}(x,t) = 10 \, x_{1} \, x_{2} \, (x_{1} - 1) \, (x_{2} - 1) \left[\frac{\sin(2 \, x_{1} \, t)}{\exp(t/2)} + \frac{\cos(x_{2} \, t)}{\exp(t/4)} + \frac{\sin(x_{1} \, x_{2} \, t)}{\exp(t)}\right].
\end{gather*}
Here, we have $ a(x) \equiv \nu $. We modify the weak formulation by partially integrating the terms with derivatives with respect to $ x_{1} $ and $ x_{2} $. This leads to the following new term for the SGA:
\begin{gather*}
    \int_{\Omega} \left(\partial_{x_{1}} u(x,t)^{2} + \partial_{x_{2}} u(x,t)^{2}\right) \phi(x) \dx{x} = -\int_{\Omega} u(x,t)^{2} \left(\partial_{x_{1}} \phi(x) + \partial_{x_{2}} \phi(x)\right) \dx{x}.
\end{gather*}
For the finite element discretization, this results in a term quite similar to $ \vec{\ell}(c,u;\mu) $:
\begin{gather*}
    \left[\vec{m}(d,\vec{u};t)\right]_{i} \mdef \int_{\Omega} d(x,u_{h}(x),\nabla u_{h}(x);t) \left( \partial_{x_{1}} \phi_{i}(x) + \partial_{x_{2}} \phi_{i}(x) \right) \dx{x},
\end{gather*}
where $ d(x,u,\nabla u;t) = u(x,t)^{2} $ in this example. However, since $ d $ is quadratic, we also consider a SGA using a multi-linear form as discussed in \cite{Tolle2020}. This leads to a simple reduction of the problem without needing complexity reduction. The SGA and multi-linear standard Galerkin approach (ML-SGA) solve the following semi-discrete problems:
\begin{alignat*}{2}
    \mat{E} \, \partial_{t} \vec{u}(t) + \nu \, \mat{K} \, \vec{u}(t) &- \frac{1}{2} \, \vec{m}(\vec{u}(t)) &~=~ \vec{q}(t), \\
    \mat{E} \, \partial_{t} \vec{u}(t) + \nu \, \mat{K} \, \vec{u}(t) &- \frac{1}{2} \, \tens{M} : (\vec{u}(t) \otimes \vec{u}(t)) &~=~ \vec{q}(t),
\end{alignat*}
respectively, where $ \mat{E} $ and $ \mat{K} $ denote the standard mass and stiffness matrices and the third-order tensor $ \tens{M} $ is defined through 
\begin{gather*}
    [\tens{M}]_{i,j,k} = \int_{\Omega} \phi_{k}(x) \, \phi_{j}(x) \, \left(\partial_{x_{1}} \phi_{i}(x) + \partial_{x_{2}} \phi_{i}(x)\right) \dx{x}.
\end{gather*}
The ML-SGA has the same benefits as the EGFEM, i.e., the cost-intensive numerical integration only has to be performed once. In contrast, the EGFEM solves the following semi-discrete problem:
\begin{gather*}
    \mat{E} \, \partial_{t} \vec{u}(t) + \nu \, \mat{K} \, \vec{u}(t) - \frac{1}{2} \, \mat{M}^{c} \, \vec{d}(\vec{u}(t)) = \vec{q}(t),
\end{gather*}
where $ \vec{d} $ is defined in a similar manner to $ \vec{a} $ and $ \vec{c} $ in \eqref{eqn:egfem}. In addition to the reduced terms in \eqref{eqn:reduced_sga} and \eqref{eqn:reduced_egfem}, we also introduce the reduced mass and stiffness matrices as well as the reduced tensor for the ML-SGA, which read:
\begin{align*}
    \mat{E}_{\mathsf{r}} &= \mat{V}_{u}^{\transpose} \, \mat{E} \, \mat{V}_{u}, &
    \mat{K}_{\mathsf{r}} &= \mat{V}_{u}^{\transpose} \, \mat{K} \, \mat{V}_{u}, &
    \tens{M}_{\mathsf{r}} = \tens{M} \times_{1} \mat{V}_{u}^{\transpose} \times_{2} \mat{V}_{u}^{\transpose} \times_{3} \mat{V}_{u}^{\transpose}.
\end{align*}
The advantage of the ML-SGA is that the nonlinearity no longer depends on the full dimension, i.e., $ \tens{M}_{\mathsf{r}} \in \rr^{n_{u} \times n_{u} \times n_{u}} $, which means that the DEIM is not needed for this formulation. However, the ML-SGA is only possible for polynomial nonlinearities and only practical for low-degree polynomials.

\begin{figure} \centering
\includegraphics[width=0.8\linewidth]{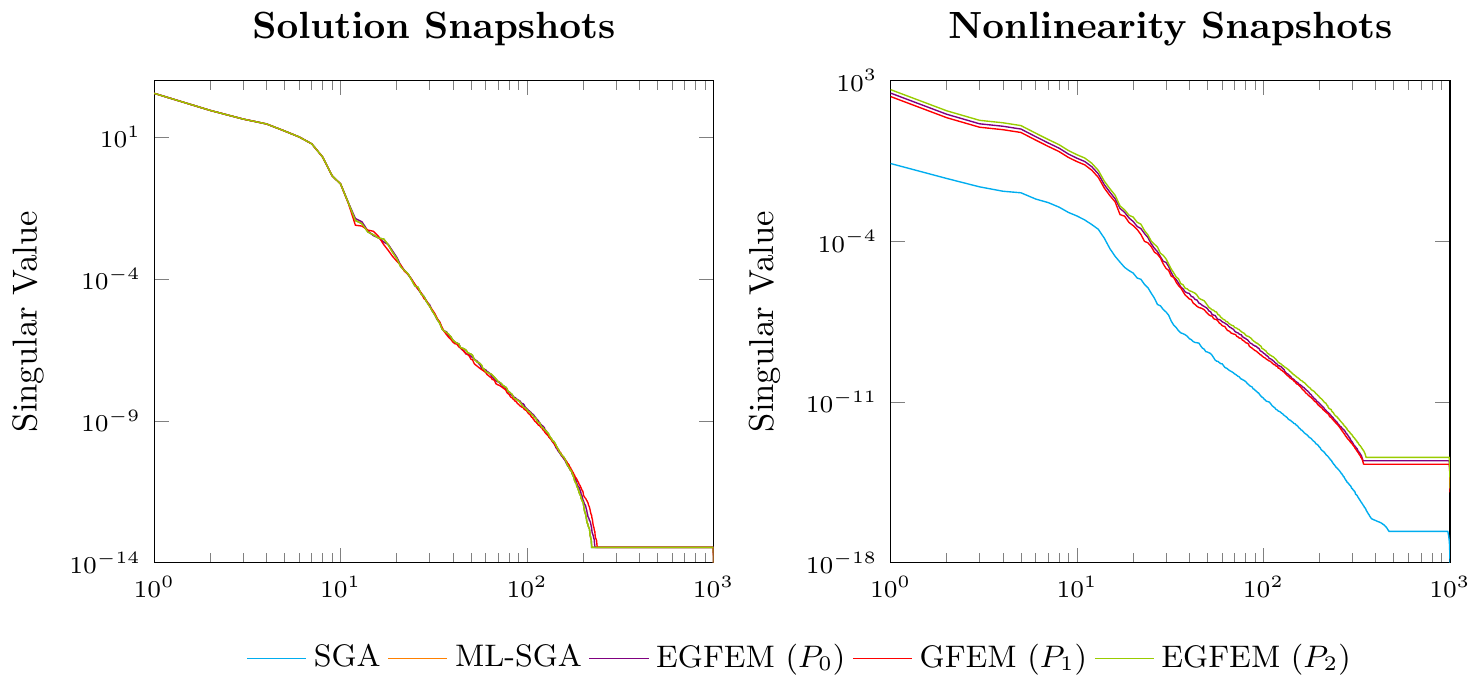}
\caption{Singular values associated to the snapshot matrix containing the full-order solution (\textsl{left}) and nonlinearity (\textsl{right}).}
\label{fig:burgersModes}
\end{figure}

\begin{figure} \centering
\includegraphics[width=0.8\linewidth]{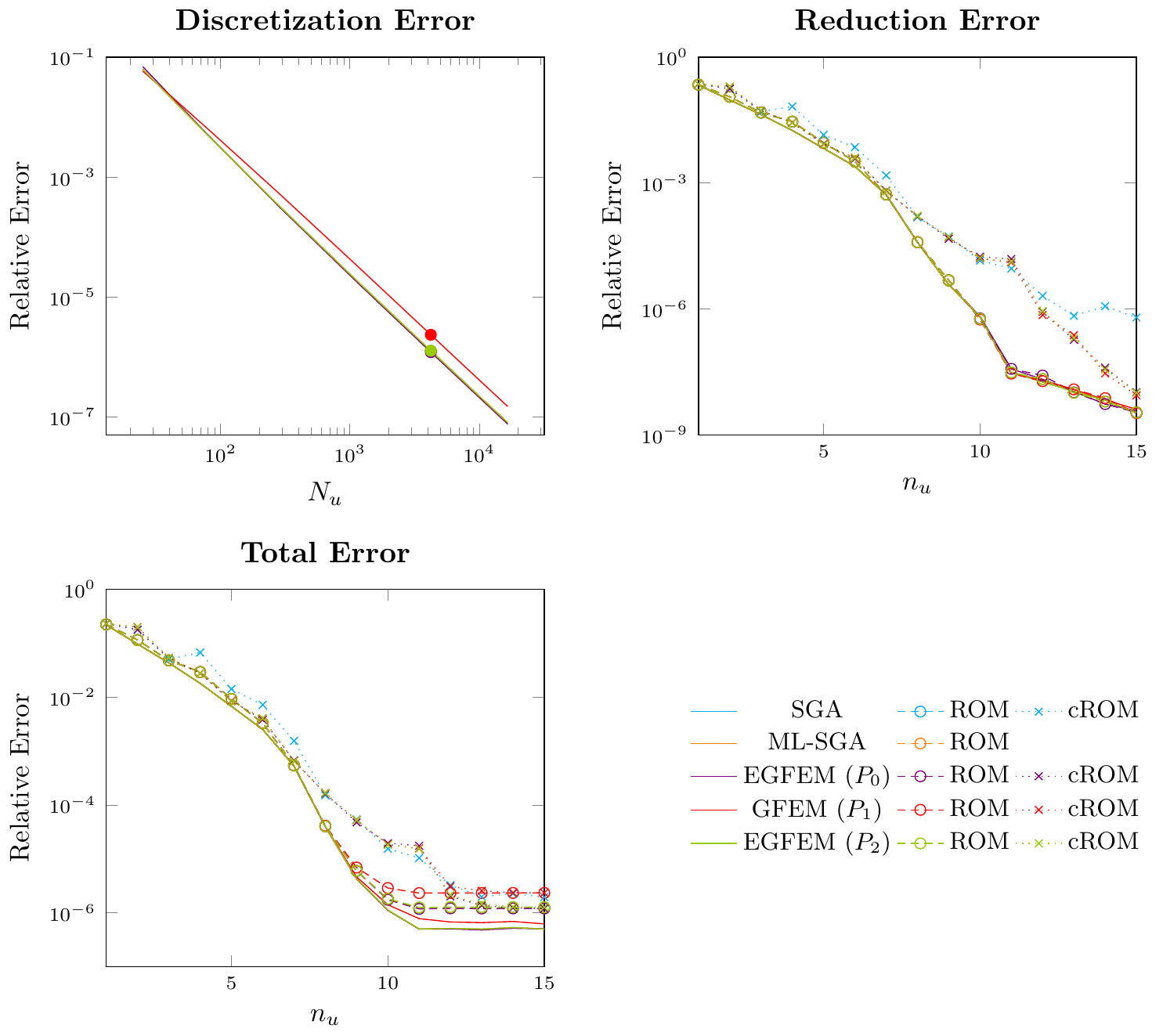}
\caption{Discretization error with respect to mesh refinement (\textsl{top, left}), where the colored marker denotes the mesh used for reduction. The reduction error (\textsl{top, right}) compares the reduced order solutions to the full order solutions, while the total error (\textsl{bottom, left}) compares the reduced order solutions to the analytical solution and includes the discretization error.}
\label{fig:burgersErrors}
\end{figure}

We begin by examining the singular values of the snapshot matrices. The snapshots are generated by solving the FOMs and evaluating the solutions on an equidistant time grid with step size $ \delta t = 10^{-2} $. Figure \ref{fig:burgersModes} displays the singular values associated with the respective snapshot matrix. We see that the singular values of the various formulations are qualitatively similar. For this example, we choose the number of DEIM modes to be equal to the number of POD modes. Because the GFEM loses accuracy with respect to the energy norm in comparison to the SGA, see \cite{Christie1981}, we use the energy norm to compute the errors. Relative errors are shown in Figure~\ref{fig:burgersErrors}. The discretization error compares the analytical solution $ \hat{u} $ to the solution of the FOM for increasingly finer triangulations. The reduction error compares the FOMs with $ N_{u} = 4,225 $ to the ROMs. The solid lines without markers show the projection error, i.e., the error arising from restricting and prolongating the full-order solution. Finally, the total error compares the ROMs to the analytical solution. From the reduction error in Figure~\ref{fig:burgersErrors}, we see that the SGA-cROMs perform worse than the other cROMs as the size of the reduced models increases. However, when considering the total error, which includes the discretization error, we see that the differences become less obvious. We do note that the EGFEM models perform slightly better than the GFEM models, which shows that our method can be better than the existing approaches using GFEM. To summarize, this problem is most easily reduced using the ML-SGA, which is applicable for low-degree polynomial nonlinearities. However, the EGFEM using quadratic ($ P_{2} $) elements is analytically equivalent to the SGA, while also being simpler to reduce and much more efficient in terms of the computational costs.

\begin{remark}[Instability]
When using DEIM, we have found that certain choices of POD and DEIM sizes may lead to unstable approximations. In particular, the cROM may be instable, when the DEIM dimension is smaller than the POD dimension. This phenomenon has also been reported for other models, see, for example, \cite{Antil2014} and \cite{Sipp2020}.
\end{remark}

\subsection{Minimal Surface Equation}

Finally, we consider an example for which the GFEM is not applicable. The minimal surface equation with a source term \cite{Verfuerth1994} is given through
\begin{align*}
    a(\nabla u) &= \left(1 + \| \nabla u \|^{2}\right)^{-\frac{1}{2}}, &
    c(x) &\equiv 0, & 
    q(x;\mu) &= 2 \left(e^{\|x\|^{2}} - e^{-1}\right) e^{\mu_{1} \, x_{1} + \mu_{2} \, x_{2}}
\end{align*}
with  $ \mu = (\mu_{1},\mu_{2}) \in \mathcal{P} = [0,1]^{2} $. The domain $ \Omega $ is a polygonal approximation of the unit disk. The nonlinear function $ a $ is an example of a nonlinearity, that cannot be approximated using the GFEM. Since $ V_{h} \subset H^{1}_{0} $, the gradient $ \nabla u_{h} $ of $ u_{h} \in V_{h} $ is not well-defined on the nodes. This means that the interpolation of $ a(\nabla u_{h}) $ onto $ V_{h} $ is not possible and as such any reduction methods based on the GFEM are not applicable. Therefore, we only consider EGFEM with a piecewise constant ($ P_{0} $) approximation of the nonlinearity. Interestingly, this is an exact reformulation of the nonlinearity, since $ \nabla u_{h} $ is piecewise constant.

\begin{figure} \centering
\includegraphics[width=0.8\linewidth]{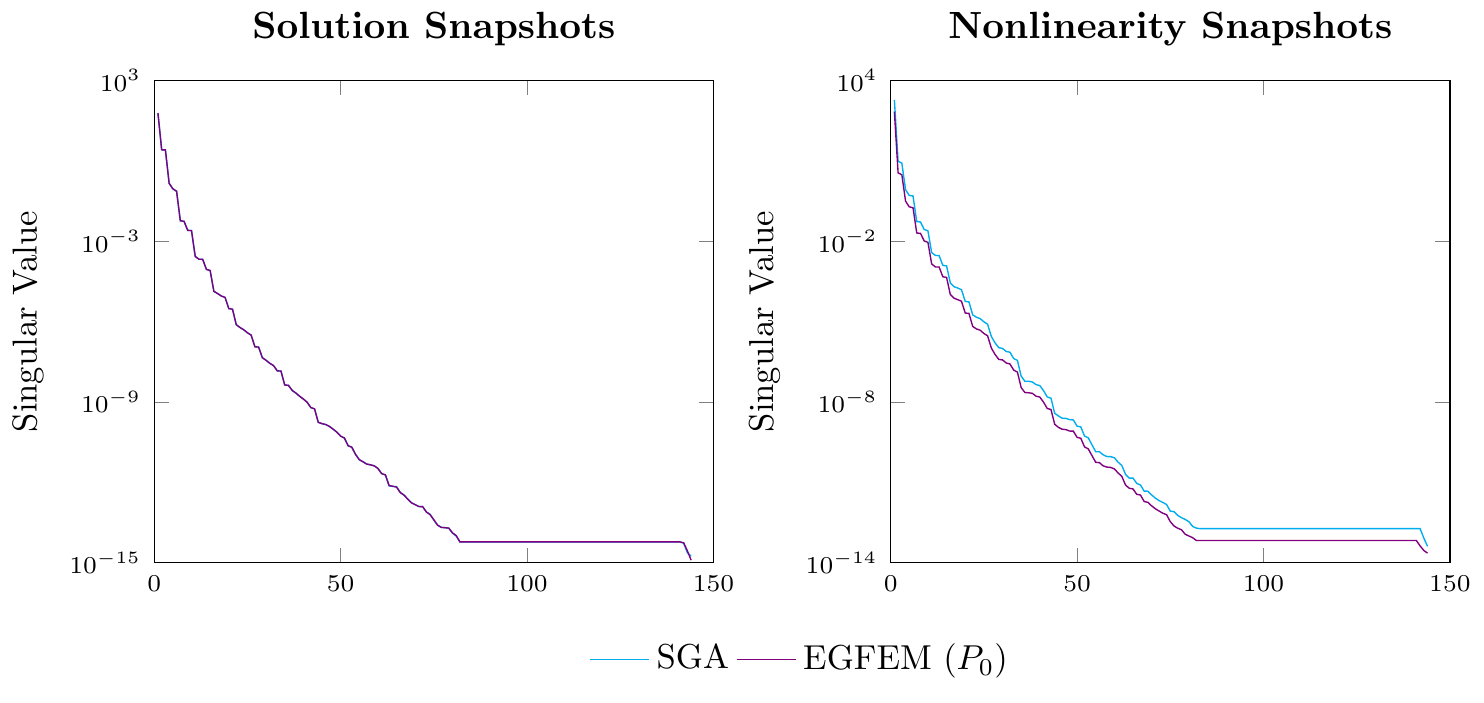}
\caption{Singular values associated to the snapshot matrices for the solution (\textsl{left}) and the nonlinearity (\textsl{right}).}
\label{fig:minsurfaceModes}
\end{figure}

\begin{figure} \centering
\includegraphics[width=0.8\linewidth]{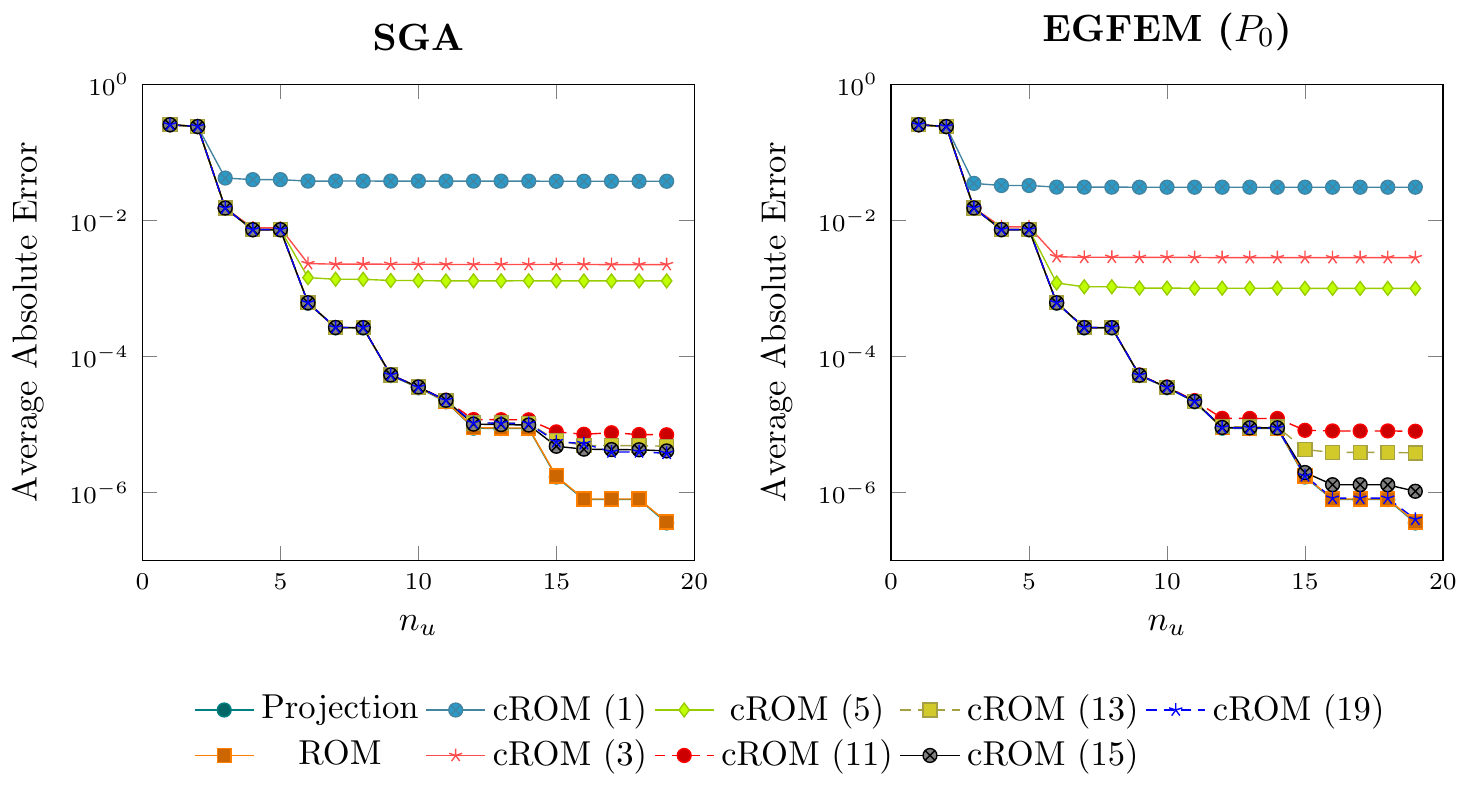}
\caption{Average absolute error over all samples between the reduced order and full order models for the SGA (\textsl{left}) and EGFEM (\textsl{right}). The errors are computed using the Euclidean (point-wise) norm.}
\label{fig:minsurfaceError}
\end{figure}

\begin{figure} \centering
\includegraphics[width=0.8\linewidth]{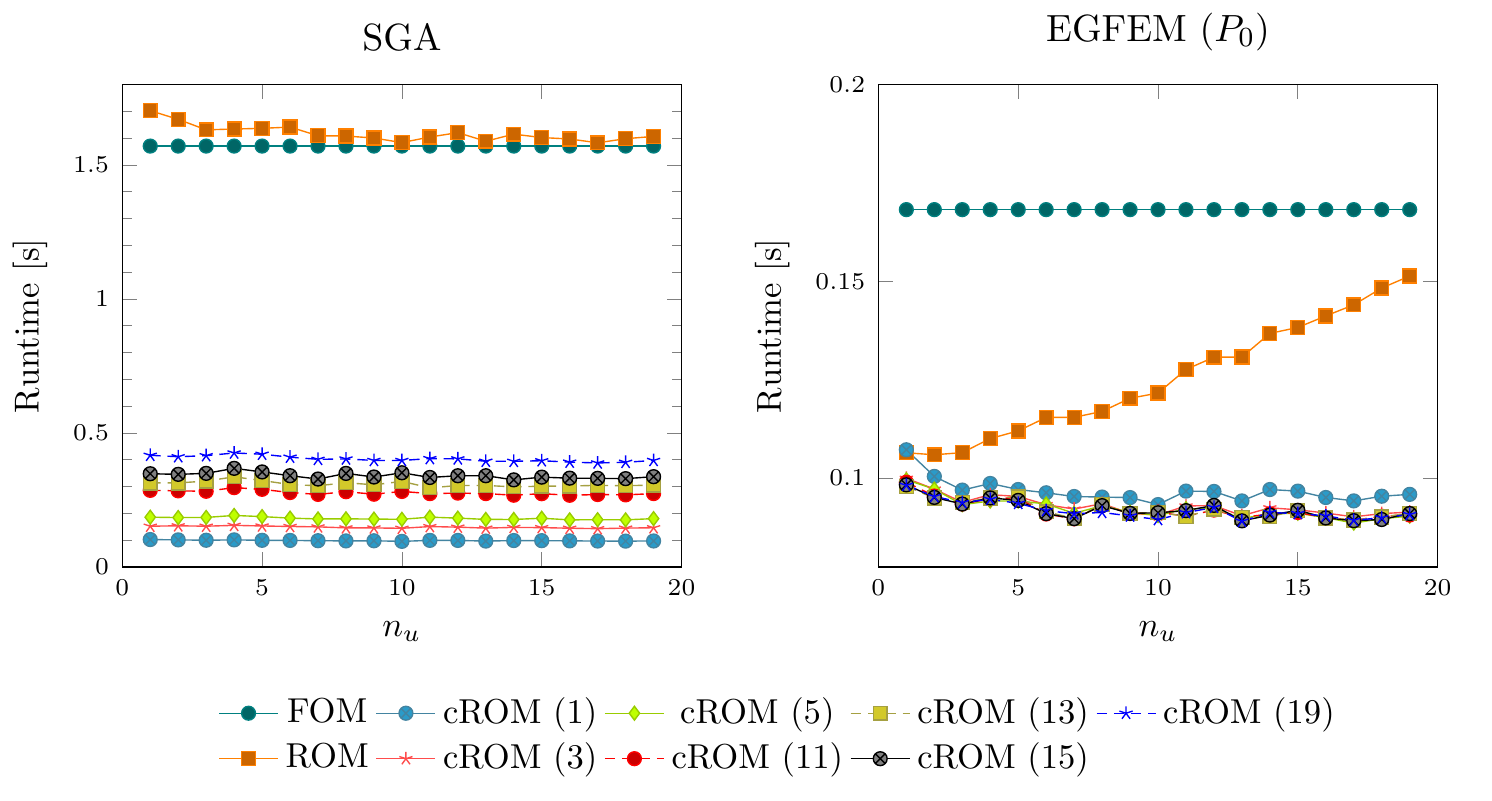}
\caption{Average time per iteration associated with the respective system.}
\label{fig:minsurfaceTimes}
\end{figure}

Similar to the first example, we sample the parameter space $ \mathcal{P} $ uniformly at 144 points, while the errors are computed on a uniform sampling with 225 points, so that the ROMs and cROMs are evaluated on parameter settings not included in the training set. The singular values for the solution and nonlinear snapshot matrices are shown in Figure~\ref{fig:minsurfaceModes}. From the singular values, we expect good approximations using few modes. The average absolute error over all samples is shown in Figure~\ref{fig:minsurfaceError}. The errors of both methods are comparable, although the EGFEM shows smaller errors as the size of the cROM increases. We also compare the average time per iteration in Figure~\ref{fig:minsurfaceTimes}. Once again, we see that the ROM without complexity reduction takes as long as the FOM for the SGA. The cROMs show a noticeable reduction in the computational overhead associated with evaluating the nonlinearity. When comparing the two approaches, we see that the EGFEM is noticeably faster than the SGA-counterparts. In summary, the EGFEM leads to faster models that are easier to reduce and prove to be more accurate than the standard approach.

\section{Conclusion} \label{sec:conclusion}

In this work, we combine our recently introduced EGFEM with nonlinear reduction techniques. The presented results use POD and DEIM. However, the choice is illustrative and not fixed. Note that our approach is also applicable with other model order and complexity reduction methods, such as the reduced basis method in place of POD and the empirical interpolation method instead of the DEIM. Due to the nature of the EGFEM, we are able to effectively reduce quasilinear problems as well as handle discontinuous nonlinearities that could not be treated using the methods presented in \cite{Dickinson2010} and \cite{Wang2015}. At the same time, the same approaches are embedded in our framework, since EGFEM is a generalization of the GFEM, which results for the special choice $ W_{h} = V_{h} $. The superior performance of our method in comparison to existing methods has been numerically demonstrated for different applications in Section~\ref{sec:results}. Furthermore, it is less intrusive in terms of an efficient implementation of the complexity reduction in comparison to standard reduction approaches.

\bibliographystyle{siam}
\bibliography{references.bib}

\end{document}